\documentclass[11pt,reqno,a4paper]{amsart}

\usepackage{amsmath,amssymb,amsthm,amsaddr}
\usepackage{enumitem}
\usepackage[margin=3cm,top=4cm,bottom=4cm,footskip=1cm,headsep=2cm]{geometry}
\usepackage[utf8]{inputenc}
\usepackage{amsthm}
\usepackage[british]{babel}

\usepackage{accents}
\newcommand{\ubar}[1]{\underaccent{\bar}{#1}}

\author{H. Egger \and J. Giesselmann}
\address{Department of Mathematics, TU Darmstadt, Germany}
\email{herbert.egger@tu-darmstadt.de}
\email{jan.giesselmann@tu-darmstadt.de}

\title[Relative energy estimates for gas transport]{Stability and asymptotic analysis for instationary gas transport via relative energy estimates}

\newtheorem{theorem}{Theorem}

\newtheorem{lemma}[theorem]{Lemma}
\newtheorem{notation}[theorem]{Notation}
\newtheorem{corollary}[theorem]{Corollary}

\theoremstyle{definition}
\newtheorem{remark}[theorem]{Remark}
\newtheorem{assumption}[theorem]{Assumption}

\def\ddt{\frac{d}{dt}}
\def\ddtau{\frac{d}{d\tau}}
\def\dt{\partial_t}
\def\dtau{\partial_{\tau}}
\def\dx{\partial_x}
\def\E{\mathcal{E}}
\def\H{\mathcal{H}}

\def\R{\mathcal{R}}
\def\J{\mathcal{J}}
\def\B{\mathcal{B}}
\def\Bd{\mathcal{B}_\partial}
\def\G{\mathcal{G}}
\def\V{\mathcal{V}}
\def\C{\mathcal{C}}
\def\D{\mathcal{D}}
\def\P{\mathcal{P}}
\def\Pd{\mathcal{P}_\partial}

\def\z{z}

\def\Hu{\H'}
\def\Huu{\H''}

\def\DD{\mathbb{D}}
\def\VV{\mathbb{V}}
\def\WW{\mathbb{W}}

\def\RR{\mathbb{R}}

\def\la{\langle}
\def\ra{\rangle}

\def\eps{\varepsilon}

\newcommand{\bu}{{\boldsymbol u}}
\newcommand{\bv}{{\boldsymbol v}}
\newcommand{\bw}{{\boldsymbol w}}
\newcommand{\bz}{{\boldsymbol z}}

\newcommand{\buh}{\widehat\bu}
\newcommand{\be}{{\boldsymbol r}}
\newcommand{\beh}{\widehat\be}

\newcommand{\bzd}{\bz_{\boldsymbol \partial}}

\newcommand{\hu}{\widehat{\boldsymbol u}}

\usepackage{xcolor}

\begin{document}

\textit{\vskip-3em \noindent Draft: \today}

\begin{abstract}
We consider the transport of gas in long pipes and pipeline networks for which the dynamics are dominated by friction at the pipe walls. The governing equations can be formulated as an abstract dissipative Hamiltonian system which allows us to derive perturbation bounds by means of relative energy estimates. As particular consequences, we obtain stability with respect to initial conditions and model parameters and quantitative estimates in the high friction limit. Our results are established in detail for the flow in a single pipe and through the energy-based modelling they naturally generalize also to pipe networks. 
\end{abstract}

\maketitle

\begin{quote}
{\footnotesize 

}
\end{quote}
\bigskip 

\section{Introduction} \label{sec:1}

Gas transport through long pipes is usually dominated by friction at the pipe walls. On the practically relevant time scales $t=O(1/\eps)$, the governing balance equations for mass and momentum can be phrased in rescaled form as 
\begin{align}
a \dtau \rho + \dx m &= 0  \label{eq:sys1}\\
\eps^2 \dtau w + \dx h + 
\gamma |w| w 
&= 0. \label{eq:sys2}
\end{align}
Here $a$ is the cross section of the pipe, $\rho$ is the density of the gas, and $w$, $\tau$, $\gamma>0$ denote the rescaled velocity, time, and friction coefficient, respectively. Furthermore 
\begin{align}
m = a \rho w 
\qquad \text{and} \qquad 
h = \frac{\eps^2 w^2}{2} + P'(\rho) + g z \label{eq:sys3}
\end{align}
are the rescaled mass flow rate and total specific enthalpy with $P(\rho)$ denoting the pressure potential, $g$ the gravity constant, and $z$ the elevation of the pipe. 
For convenience of the reader, a detailed derivation of the above equations, which are assumed to hold for $0 < x < \ell$ and for time $\tau \ge 0$, is given in the appendix. 
In rescaled variables, the physical energy of the system is given by 
\begin{align}
\H(\rho,w) &= \int_0^\ell a \left( \eps^2 \frac{\rho  w^2}{2} + P(\rho) + g z \rho \right) \, dx. \label{eq:sys4}
\end{align}
Let us note that the state variables $(\rho,w)$ and the co-state variables $(h,m)$ are directly linked via the derivative of this energy functional; details will be given below. 
As a direct consequence of the particular problem structure, sufficiently smooth solutions of the system \eqref{eq:sys1}--\eqref{eq:sys3} can be shown to satisfy the following energy-dissipation inequality
\begin{align}
\ddtau \H(\rho,w)  + \D(\rho,w) &= -m h \big|_0^\ell \label{eq:sys5}
\end{align}
with dissipation functional $\D(\rho,w) = \int_0^\ell \gamma \rho |w|^3 \, dx \ge 0$. 
The free energy of the system thus changes only due to friction at the pipe walls and energy transfer across the boundary. 

To fully describe the evolution, the problem has to be complemented by appropriate initial and boundary conditions. For instance, one may prescribe  
\begin{align} \label{eq:sys6}
h = h_\partial \qquad \text{on } \{0,\ell\};
\end{align}
alternative conditions will be discussed later. We will however mostly study properties of general solutions of \eqref{eq:sys1}--\eqref{eq:sys3} without restriction of the boundary values.

\subsection*{Scope and summary of main results.}

In this paper, we investigate the stability of solutions to the nonlinear system \eqref{eq:sys1}--\eqref{eq:sys3} with respect to perturbations of the initial and boundary values as well as the model parameters $\eps$ and $\gamma$. 
Our analysis is done under the following structural assumptions:
\begin{itemize}
\item[(A1)] The pressure potential $P: \RR_+ \to \RR$ is smooth and strictly convex and there exist positive constants $\ubar \rho,\bar \rho, \bar w$ and $\bar \eps$ such that 
        \begin{align} \label{eq:ass1}
        \rho P''(\rho) \ge 4 \bar \eps^2 |\bar w|^2 \qquad \forall  \ubar \rho \le \rho \le \bar \rho.
    \end{align}
    Moreover, 
    $0<\ubar a \le a \le \bar a$ and $|gz| \le \bar g\bar z$ for appropriate constants $\ubar a,\bar a$, and $\bar g \bar z$. 
\item[(A2)] Sufficiently smooth solutions $(\rho,w)$ and $(\hat \rho,\hat w)$ to the system \eqref{eq:sys1}--\eqref{eq:sys5} exist for parameters $0 \le \eps, \hat \eps \le \bar \eps$ and $0 < \ubar \gamma \le \gamma,\hat \gamma \le \bar \gamma$ with $\ubar \gamma,\bar \gamma$ constant, which satisfy
\begin{align} \label{eq:ass2}
    \ubar{\rho} \le \rho,\hat \rho \le \bar \rho
    \qquad \text{and} \qquad 
 -  \bar w \le w,\hat w \le\bar w.
\end{align}
\end{itemize}
Conditions (A1) and (A2) imply uniform bounds for $P$ and its derivatives and ensure that the flow is subsonic; we refer to \cite{Dafermos2016} for details. 
Under these conditions, we will show the following stability result:  
Let $(\rho,w)$ and $(\hat \rho,\hat w)$ be sufficiently regular solutions of \eqref{eq:sys1}--\eqref{eq:sys3} with parameters $\eps,\gamma$ and $\hat \eps,\hat \gamma$, and boundary values $\hat h_\partial$ and $h_\partial$ as described in \eqref{eq:sys6}. Then
\begin{align*}
    \|\rho(\tau) - \hat \rho(\tau)\|_{L^2}^2 &+ \eps^2 \|w(\tau) - \hat w(\tau)\|^2_{L^2} + \int_0^\tau \|w(s)-\hat w(s)\|^3_{L^3} ds \\
    &\le \hat C e^{\hat c \tau} \Big(\|\rho(0) - \hat \rho(0)\|_{L^2}^2  + \|\eps w(0) - \eps \hat w(0)\|_{L^2}^2  \\
    & \qquad \qquad \qquad +  |\gamma - \hat \gamma|^{3/2} + |\eps^2 -\hat \eps^2| + \int_0^\tau |h_\partial(s) - \hat h_\partial(s)| \, ds \Big),
\end{align*}
with constants $\hat c,\hat C$ depending only on the bounds in assumptions (A1)--(A2) and on bounds for time derivatives of $\hat \rho$ and $\hat w$.   
Note that the energy \eqref{eq:sys4} and dissipation functional \eqref{eq:sys5}, and as a consequence also the co-state variables \eqref{eq:sys3} depend explicitly on the parameters $\eps$ and $\gamma$, so their definition changes for the perturbed problem.
A precise statement of our results and of the regularity assumption on the solutions is given in Section~\ref{sec:4}, where we also discuss various generalizations including perturbations in other model parameters and the choice of different boundary conditions. 
Let us mention two immediate consequences of the above stability estimate, namely
\begin{itemize}
\item uniqueness of regular subsonic bounded state solutions for specified initial and boundary values,  and their stability with respect to perturbations in these problem data as well as the model parameters; 
\item convergence of solutions to those of the parabolic limit problem which results by formally setting $\eps=0$ in equations \eqref{eq:sys1}--\eqref{eq:sys6}. 
\end{itemize}
Existence and uniqueness of solutions to the parabolic limit problem can be proven rigorously by variational arguments; see \cite{Bamberger79,Raviart70} or \cite{SchoebelKroehn2020}. 
This parabolic problem also serves as the basis for simulation codes utilized in the gas network community \cite{BrouwerGasserHerty2011,Osiadacz87}, and the stability estimate above allows us to obtain a quantitative justification for the use of this model. 
Due to the energy-based modelling, the results can be generalized almost verbatim to gas networks by utilizing appropriate coupling conditions; details will be discussed in Section~\ref{sec:5}.

\subsection*{Main tools.}

The proof of our main result is based on the observation that \eqref{eq:sys1}--\eqref{eq:sys3} can be written as an abstract port-Hamiltonian system \cite{VanDerSchaftMaschke2002}
\begin{align} \label{eq:dhs}
    \C \dtau \bu + (\J+\R(\bu)) \bz(\bu) &= \Bd \bz(\bu), 
\end{align}
with $\bu=(\rho,w)$ and $\bz(\bu)=(h,m)$ denoting the state and co-state variables, which are linked via the energy functional $\H(\bu)=\H(\rho,w)$; 
moreover, $\C$, $\J$, $\R(\bu)$, and $\Bd$, are appropriate operators, the last one incorporating the boundary conditions; see Section~\ref{sec:2}.

Any smooth function $\buh=(\hat \rho,\hat w)$, e.g., a solution of \eqref{eq:sys1}--\eqref{eq:sys3} with perturbed model parameters, may be considered as a solution of the system
\begin{align} \label{eq:dhsp}
  \C \dtau \buh + (\J + \R(\buh)) \bz(\buh) &= \Bd \bz(\buh) + \beh,
\end{align}
with the same operators $\C$, $\J$, $\R(\cdot)$, $\Bd$, and the same energy functional $\H(\cdot)$ and state to co-state mapping $\bz(\cdot)$, up to some perturbation $\beh$.
The coincidence of the underlying Hamiltonian structure will allow us to estimate the difference between $\bu$ and $\buh$ in terms of the perturbations $\beh$ by means of relative energy estimates.
For the stability analysis on the abstract level, we require some general conditions that will be verified in detail for the gas transport problem under investigation using the assumptions (A1)--(A2).

The use of an abstract problem structure greatly simplifies the analysis and allows to generalize our results in various directions. We briefly discuss the incorporation of perturbations in other model parameters and the extension to gas networks.

\subsection*{Review of related work}

Relative entropy or energy techniques have been used intensively for the existence, stability, and discretization error analysis of time dependent partial differential equations. 
We refer to \cite{Juengel2016} for a recent summary of corresponding results for parabolic evolution problems. 
In the present paper, we are interested in hyperbolic problems, where the use of relative entropy arguments goes back to the seminal works of DiPerna \cite{DiPerna79} and Dafermos \cite{Dafermos79}; also see \cite{Dafermos2016} for an introduction to the field.
Typical aspects that are addressed are: convergence to steady states, stable dependence of solutions on initial data and parameters, and asymptotic limits.
Examples for the latter include the low Mach limit of Euler and Navier-Stokes equations, which are investigated in \cite{FeireislNovotny2017} for instance. 
Long time convergence of solutions to damped Euler equations to Barenblatt solutions were investigated by Huang and coworkers in a series of papers \cite{GengHuang2019}.

One particular aspect that we want to address in the present study are parabolic limits of quasilinear hyperbolic equations; see  \cite{MarcatiMilani1990,JuncaRascle2002,LattanzioTzavaras2013,LattanzioTzavaras2017} for some exemplary results in this direction. 
The latter reference as well as \cite{CarrilloPengWroblewskaKaminska2020,GiesselmannLattanzioTzavaras2017} strongly rely on the underlying dissipative Hamiltonian structure of many equations in fluid mechanics, which will also play a major role in our analysis below. 
The previous papers, however, use formulations in conservative variables, which allows to deal with one solution being a weak solution in the classical sense of hyperbolic conservation laws \cite{Dafermos2016}. 
Parabolic limits of hyperbolic $2\times2$ systems have also been studied in \cite{LinCoulombel2013,XuKawashima2014} using compensated compactness arguments \cite{HuangMarcatiPan2005, HuangPanWang2011}. This has the advantage that no additional regularity of solutions to the parabolic limit problem is required, but no quantitative information about the speed of convergence is obtained. 
Spectral estimates for the linear part of the problem are used in \cite{DuanLiuZhu2015} to derive convergence in the parabolic limit. 
Let us note that most of the above works consider only linear friction laws and unbounded domains or periodic boundary conditions. 

In contrast to work mentioned previously, our study is based on a formulation in primitive variables, which requires both solutions to have some minimal smoothness. On the other hand, this formulation allows us to incorporate boundary conditions more naturally and to extend our results to networks in a straight-forward manner using appropriate coupling conditions at network junctions \cite{Reigstad2014,Egger2018} that guarantee energy conservation or dissipation at network junctions.
Similar formulations for compressible flow were also considered in the context of port-Hamiltonian systems; see \cite{VanDerSchaftMaschke2002,VanDerSchaftJeltsema2014} for the models and \cite{Cardoso2019,LiljegrenSailer2020} for corresponding discretization strategies. 
Other systems, that fit into the general framework that we develop in this paper include the Euler-Korteweg system, the system of quantum hydrodynamics and the Euler-Poisson equations. An overview about corresponding results can be found in \cite{GiesselmannLattanzioTzavaras2017}.

%

\subsection*{Outline of the manuscript.}

The remainder of the manuscript is organized as follows:
Section~\ref{sec:2} is concerned with the perturbation analysis for the abstract system \eqref{eq:dhs}. 
In Section~\ref{sec:3}, we then briefly review the derivation of the model equations \eqref{eq:sys1}--\eqref{eq:sys3} starting from the usual balance equations of gas dynamics, and we show that they fit into the abstract form \eqref{eq:dhs}.
In Section~\ref{sec:4}, we verify the assumptions required for our abstract analysis for the gas transport problem under investigation, which allows us to state and prove our main results. 
Their generalization to gas networks is discussed in Section~\ref{sec:5}.
%


\section{An abstract stability estimate} \label{sec:2}


In this section, we consider abstract evolution problems of the form 
\begin{align}  \label{eq:abs1}
\C \dtau \bu + (\J + \R(\bu)) \, \bz(\bu) = \Bd \bz(\bu),  \\
\bz(\bu) = \C^{-1} \H'(\bu), \label{eq:abs2}
\end{align}
with state and co-state variables $\bu$ and $\bz(\bu)$ that are directly connected via the derivative $\H'(\bu)$ of an associated energy functional $\H(\bu)$. 
After briefly introducing a reasonable functional analytic setting, we derive stability estimates for solutions.

\subsection{Notation and basic assumptions}

Let $\WW \subset \VV$ be real Hilbert spaces, with $\VV' \subset \WW'$ denoting the corresponding dual spaces. We use $\la \cdot,\cdot \ra$ to denote the duality product on $\VV' \times \VV$ and $\WW' \times \WW$.
Furthermore, let $\H : \DD \subset \VV \to \RR$ be a smooth and strictly convex energy functional, with $\DD \subset \VV$ denoting some appropriate, e.g., convex and closed, subset. 
We consider abstract evolution problems of the form 
\eqref{eq:abs1}--\eqref{eq:abs2},
with operators satisfying the following conditions. 
\begin{assumption} \label{ass:1}
$\C:\VV \to \VV'$ is linear, bounded, self-adjoint, and elliptic on $\VV$, i.e., 
    \begin{align}
       \la \C \bu, \bv \ra = \la \C \bv, \bu\ra \qquad &\forall \bu,\bv \in \VV, \label{eq:assC1}\\
       c \|\bv\|_\VV^2 \le \la \C \bv, \bv\ra \le C \|\bv\|_{\VV}^2 \qquad &\forall \bv \in \VV, \label{eq:assC2}
    \end{align}
    with positive constants $c,C$ independent of $\bv$.
For any $\bu \in \DD$ the operator $\R(\bu) : \WW \to \WW'$ is linear, bounded, self-adjoint, and non-negative, i.e.,
\begin{alignat}{5}
    \la \R(\bu) \bw, \bz \ra &= \la \R(\bu) \bz, \bw\ra \qquad &&\forall \bw,\bz \in \WW, \label{eq:assR1}\\
    \la \R(\bu) \bw, \bw\ra &\ge 0  \qquad &&\forall \bw \in \WW.\label{eq:assR2}
\end{alignat}
Furthermore, $\J : \WW \to \WW'$ is linear and 
anti-symmetric, i.e., 
  \begin{align} 
     \la \J \bw, \bz \ra = - \la \J \bz, \bw \ra \qquad \forall \bw,\bz \in \WW, \label{eq:assJ} 
  \end{align}
and finally, the operator $\Bd: \WW \to \WW'$ is linear and bounded. 
\end{assumption}
%
From conditions \eqref{eq:assC1}--\eqref{eq:assC2}, one can immediately see that 
\begin{align} \label{eq:scalprod}
\la \bu, \bv \ra_\C:=\la \C \bv, \bu \ra
\end{align}
defines a scalar product on $\VV$ and the associated norm $\|\bv\|_\C^2 = \la \bv,\bv\ra_\C =  \la \C \bv, \bv \ra$ is equivalent to the standard norm $\|\cdot\|_\VV$. 
The expression $\bz(\bu) = \C^{-1} \Hu(\bu) = \operatorname{grad}_\C \H(\bu)$ then denotes the gradient of the functional $\H$ at $\bu$ with respect to this scalar product. We further introduce the symbol 
\begin{align} \label{eq:hess}
\G(\bu) = \C^{-1} \Huu(\bu)
\end{align}
for the Hessian operator $\G(\bu) : \VV \to \VV'$, 
and note that 
\begin{align} \label{eq:hess2}
    \la \G(\bu) \bv, \bw \ra_\C = \la \H''(\bu) \bv, \bw \ra = \la \H''(\bu) \bw, \bv\ra = \la \G(\bu) \bw, \bv\ra_\C,
\end{align}
i.e., the Hessian is 
symmetric with respect to the scalar product induced by $\C$.
%
%
\begin{notation}
By a classical solution of \eqref{eq:abs1}--\eqref{eq:abs2} on $[0,T]$, we mean a function 
$$ 
\bu \in C^1([0,T];\VV) \cap 
C^0([0,T];\DD)  \quad \text{with} \quad  
\bz( \bu) \in C^0([0,T];\WW)
$$ 
that satisfies \eqref{eq:abs1}--\eqref{eq:abs2} for all $0 \le \tau \le T$ in the sense of $\WW'$.
\end{notation}

\subsection{Power balance}

As a direct consequence of the underlying port-Hamiltonian structure and our assumptions on the operators, we obtain the following power balance relation.
\begin{lemma} 
Let $\bu$ be a classical solution of \eqref{eq:abs1}--\eqref{eq:abs2}. Then 
\begin{align} \label{eq:powerbalance}
\ddtau \H(\bu) 
= -\D(\bu) + \la \Bd \bz(\bu), \bz(\bu)\ra.
\end{align}
with $\D(\bu) := \la \R(\bu) \bz(\bu), \bz(\bu)\ra$ denoting the dissipation functional, i.e., the total energy of the system can only change via dissipation or power flowing over the ports. 
\end{lemma}
\begin{proof}
By formal computation and equations~\eqref{eq:abs1}--\eqref{eq:abs2}, we obtain
\begin{align*}
\ddtau \H(\bu) 
&= \la \H'(\bu), \dtau \bu\ra 
 = \la \C \dtau \bu, \bz(\bu)\ra \\
&= -\la \J \bz(\bu), \bz(\bu)\ra - \la \R(\bu) \bz(\bu), \bz(\bu)\ra + \la \Bd \bz(\bu), \bz(\bu)\ra.
\end{align*}
The result then follows immediately from the assumptions on the operators.
\end{proof}



\subsection{Evolution of the relative energy}

We now study the stability of solutions to \eqref{eq:abs1}--\eqref{eq:abs2} with respect to perturbations. To this end, let $\buh$ denote a classical solution of 
\begin{align}
\C \dtau \buh + [\J+\R(\buh)] \bz(\buh) &= \Bd \bz(\bu) +  \beh, \label{eq:abs1p}\\
\z(\buh) &= \C^{-1} \H'(\buh), \label{eq:abs2p}
\end{align}
with appropriate perturbation described by the residual functional $\beh \in C^0([0,T];\WW')$. 
As a measure for the difference of $\bu$ and $\buh$, we utilize the relative energy \cite{Dafermos2016}, defined by  
\begin{align} \label{eq:relenergy}
\H(\bu|\buh) := \H(\bu) - \H(\buh) - \la \Hu(\buh), \bu-\buh \ra.
\end{align}
Using the particular problem structure and some elementary computations, we can prove the following basic identity for the temporal change of the relative entropy. 
\begin{lemma} \label{lem:abs}
Let $\bu$, $\buh$ be classical solutions of \eqref{eq:abs1}--\eqref{eq:abs2} and \eqref{eq:abs1p}--\eqref{eq:abs2p}. 
Then 
\begin{align*}
\ddtau \H(\bu|\buh)
= -\la \R(\bu) \bz(\bu) &- \R(\buh) \bz(\buh), \bz(\bu) - \bz(\buh) \ra 
  + \la \Bd (\bz(\bu) - \bz(\buh)), \bz(\bu) - \bz(\buh) \ra  \\
  & +\la \C \dtau \buh, \bz(\bu) - \bz(\buh) - \G(\buh) (\bu - \buh)\ra 
  + \la \beh, \bz(\bu) - \bz(\buh)\ra.
\end{align*}\end{lemma}
\begin{proof}
By formal differentiation of the relative energy, we obtain
\begin{align*}
\ddtau \H(\bu|\buh)
&= \la \Hu(\bu), \dtau \bu\ra - \la \Hu(\buh), \dtau \buh \ra - \la \Hu(\buh),\dtau \bu - \dtau \buh\ra - \la \Huu(\buh) \dtau \hu,  \bu - \buh\ra \\
&= \la \H'(\bu) - \H'(\buh), \dtau \bu - \dtau \buh\ra 
   + \la \H'(\bu) - \H'(\buh) - \H''(\buh) (\bu - \buh), \dtau \buh\ra \\
&= \la \C \dtau \bu - \C \dtau \buh, \bz(\bu) - \bz(\buh)\ra + \la \C \dtau \buh, \bz(\bu) - \bz(\buh) - \G(\buh) (\bu - \buh)\ra.
\end{align*}
Here we used the symmetry of $\Huu(\buh)$ in the second, and the definitions of the gradient and Hessian operators in the last step.
We now use \eqref{eq:abs1} and \eqref{eq:abs1p} to replace the time derivatives in the first term, and arrive at
\begin{align*}
\ddtau \H(\bu|\buh)
= -\la \J (\bz(\bu) &- \bz(\buh)), \bz(\bu) - \bz(\buh) \ra
   -\la \R(\bu) \bz(\bu) - \R(\buh) \bz(\buh), \bz(\bu) - \bz(\buh) \ra \\
 &  + \la \Bd (\bz(\bu) - \bz(\buh)), \bz(\bu) - \bz(\buh) \ra    + \la \beh, \bz(\bu) - \bz(\buh)\ra \\
 & +\la \C \dtau \buh, \bz(\bu) - \bz(\buh) - \G(\buh) (\bu - \buh)\ra.
\end{align*}
Due to the anti-symmetry property \eqref{eq:assJ} of the operator $\J$, the first term vanishes, and we already obtain the assertion of the lemma.
\end{proof}

\subsection{An abstract stability result} \label{sec:abs_st}

We now derive quantitative estimates for the difference of the solutions $\bu$ and $\buh$ to \eqref{eq:abs1}--\eqref{eq:abs2} and \eqref{eq:abs1p}--\eqref{eq:abs2p} with respect to perturbations in the right hand side and the initial and boundary values. 
To do so, we make some abstract assumptions that will later be verified for the gas transport problem under consideration.
\begin{assumption} \label{ass:2}
There exist constants $\hat c_0=\hat c_0(\DD)>0$, $ \hat C_0=\hat C_0(\DD)>0$ such that
\begin{align} \label{eq:norm}
\hat c_0 \|\bu - \buh\|_{\C}^2 \le \H(\bu|\buh) \le \hat C_0 \|\bu-\buh\|_{\C}^2 \qquad \text{for all } \bu, \buh \in \DD \subset \VV. \tag{C0}
\end{align}
Moreover, there exists a \emph{relative dissipation} functional $\D(\cdot|\cdot): \DD \times \DD \rightarrow [0,\infty)$ and perturbation functionals $\P(\cdot): \WW' \to \RR$ and $\Pd(\cdot): \WW \to \RR$ such that 
\begin{align}
    -\la \R(\bu) \bz(\bu) - \R(\buh) \bz(\buh), \bz(\bu)-\bz(\buh) \ra 
& \le  \hat C_1 \H(\bu|\buh) - 2 \D(\bu|\buh), \label{eq:term1} \tag{C1}\\ 
\la \C \dtau \buh, \bz(\bu) - \bz(\buh) - \G(\buh) (\bu - \buh)\ra &\le \hat C_2 \H(\bu|\buh), \label{eq:term2} \tag{C2} \\
\la \beh, \bz(\bu) - \bz(\buh)\ra
&\le \hat C_3 \H(\bu|\buh) + \D(\bu|\buh) +  \P(\beh),  \label{eq:term3} \tag{C3} \\
\la \Bd (\bz(\bu) - \bz(\buh)), \bz(\bu) - \bz(\buh) \ra  &\le \Pd(\bz(\bu) - \bz(\buh)), \label{eq:term4} \tag{C4} 
\end{align}
for classical solutions $\bu,\buh$ of \eqref{eq:abs1}--\eqref{eq:abs2} and \eqref{eq:abs1p}--\eqref{eq:abs2p} with positive constants $\hat C_i$, which may depend on the set $\DD$, the constant $\hat c_0$, $\hat C_0$, and the solution $\buh$ and its derivatives, but not on $\bu$.
\end{assumption}
Together with the relative energy identity stated in the previous lemma, these abstract conditions immediately lead to the following stability estimate.
\begin{lemma}\label{lem:abs_st}
Let Assumptions~\ref{ass:1} and \ref{ass:2} hold. Then any pair of classical solutions $\bu$ and $\buh$ to the evolution equations \eqref{eq:abs1}--\eqref{eq:abs2} and \eqref{eq:abs1p}--\eqref{eq:abs2p} satisfies
\begin{align*}
\hat c_0 \|\bu(\tau) &- \buh(\tau)\|_{\C}^2  +  \int_0^\tau e^{\hat c (\tau-\sigma)}\D (\bu|\buh) d\sigma \\
&\le \hat C_0 e^{\hat c \tau} \|\bu(0) - \buh(0)\|^2_{\C} +  \int_0^\tau e^{\hat c (\tau - \sigma)} \left[\P(\beh(\sigma)) + \Pd(\bz(\bu(\sigma)) - \bz(\buh(\sigma))) \right] \, d\sigma,
\end{align*}
with constant $\hat c = \hat C_1+\hat C_2 + \hat C_3$ and $\hat c_0,\hat C_0$ obtained from Assumption~\ref{ass:2}.
\end{lemma}
%
%
\begin{proof}
From Lemma~\ref{lem:abs} and Assumption~\ref{ass:2}, we immediately obtain
\begin{align*}
\ddtau \H(\bu|\buh) \le - \D(\bu|\buh)
+ (\hat C_1 + \hat C_2 + \hat C_3) \H(\bu|\buh) +  \P(\beh) + \Pd(\bz(\bu) - \bz(\buh)).
\end{align*}
The assertion then follows by application of the Gronwall lemma \cite[Ch.~29]{Wloka}, the definition of $\hat c$, and the estimates for the relative energy in condition \eqref{eq:norm}.
\end{proof}

\section{Application to gas networks} \label{sec:3}

We now return to the gas transport problem stated in the introduction and show that it fits into the abstract framework discussed in the previous section. In addition, we collect some auxiliary results that will be useful for the stability analysis of the next section.

\subsection{Variational formulation and canonical form}

We may multiply \eqref{eq:sys1}--\eqref{eq:sys2} by appropriate test functions $q$, $r$, integrate over the spatial domain, and use integration-by-parts in the second equation, to obtain 
\begin{align}
    (a \dtau \rho, q) + (\dx m, q) &= 0, \label{eq:var1}\\
    (\eps^2 \dtau  w, r) - (h, \dx r) + (\gamma \frac{|w|}{a\rho} m, r) &= - h r\big|_{0}^\ell, \label{eq:var2}
\end{align}
where $(a,b) := \int_0^\ell a(x) b(x) dx$ is used to abbreviate the $L^2$-scalar product. 
Note that boundary conditions \eqref{eq:sys6} for $h$ could be incorporated naturally in the last term of \eqref{eq:var2}. 
The two variational identities \eqref{eq:var1}--\eqref{eq:var2} hold for all time $t>0$ of relevance and for all smooth test functions $q$, $r$, independent of time, and they are satisfied, in particular, by all smooth solutions of \eqref{eq:sys1}--\eqref{eq:sys3}. 
In compact notation, the system \eqref{eq:var1}--\eqref{eq:var2} can be stated as
\begin{align} \label{eq:dhvar}
\la \C \dtau \bu, \bw \ra + \la \J \bz(\bu), \bw \ra + \la \R(\bu) \bz(\bu), \bw \ra &= \la \Bd \bz, \bw\ra,
\end{align}
with state variable $\bu=(\rho,w)$, co-state variable $\bz(\bu) = (h,m)$ defined by \eqref{eq:sys3}, time independent test function $\bv=(q,r)$, and operators $\C$, $\J$, $\R(\bu)$, and $\Bd$ given by 
\begin{align*}
\la \C \bu, \bv\ra &= (a \rho, q) + (\eps^2 w, r), &
\la \R(\bu) \bz, \bv \ra &= (\gamma \frac{|w|}{a\rho} m, r) \\
\la \J \bz, \bv\ra &= (\dx m,q) - (h, \dx r), &
\la \B \bzd, \bv\ra &= - h r\big|_0^\ell.
\end{align*}
From these variational characterizations, it is not difficult to see that $\J$ is anti-symmetric and that $\C$ and $\R(\bu)$ are symmetric and at least positive semi-definite, if the parameters $a$, $\eps^2$, $\gamma$ and the density $\rho$ are positive. The operator $\Bd$ associated with the boundary terms does not have a particular property, except being supported only at the boundary.

\begin{remark} 
Equation \eqref{eq:dhvar} is the variational form of an abstract evolution problem \eqref{eq:abs1} with state and co-state variables $\bu$ and $\bz(\bu)$, and energy functional $\H(\bu)=\H(\rho,w)$. Problem \eqref{eq:sys1}--\eqref{eq:sys3} thus corresponds to an abstract port-Hamiltonian system \eqref{eq:abs1}--\eqref{eq:abs2}.
\end{remark}

\subsection{Auxiliary results}

A quick inspection of the above derivations shows that the operators $\C$, $\R(\bu)$, and $\J$, and $\Bd$, can be formally identified with
\begin{align*}
    \C = \begin{pmatrix} a & 0 \\ 0 & \eps^2 \end{pmatrix}, 
\qquad 
    \R(\bu) = \begin{pmatrix} 0 & 0 \\  0 & \frac{\gamma |w|}{a \rho}\end{pmatrix},
\qquad \text{and} \qquad 
    \J -\Bd = \begin{pmatrix} 0 & \dx \\ \dx & 0 \end{pmatrix}.
\end{align*}
The latter follows rigorously by reversing the order of arguments in the derivations of the weak formulation. 
%
%
The energy of the system is here given by
\begin{align*}
\H(\bu) = \int_0^\ell a (\eps^2 \rho \frac{w^2}{2} + P(\rho) + \rho g z ) dx,
\end{align*}
and by elementary computations, we obtain the formulas
\begin{align*}
\bz(\bu) 
  = \begin{pmatrix} \eps^2 \frac{w^2}{2} + P'(\rho) \\ a \rho w \end{pmatrix} 
\qquad \text{and} \qquad 
\G(\bu) 
  = \begin{pmatrix} P''(\rho) & \eps^2 w \\ a w &  a \rho \end{pmatrix} 
\end{align*}
for the gradient $\bz(\bu) = \C^{-1} \H'(\bu)$ and Hessian $\G(\bu) = \C^{-1} \H''(\bu)$ of the energy functional.

\begin{remark}
Let us emphasize that the operators $\C$ and $\R(\cdot)$, the energy functional $\H(\cdot)$, and thus the functions $\bz(\cdot)$, $\G(\cdot)$ explicitly depend on the model parameters $\eps$ and $\gamma$. 
\end{remark}

\subsection{Functional analytic setting}\label{sec:funcana}

As a next step, we briefly discuss the choice of suitable function spaces for the gas transport problem under consideration.
We define
\begin{align} \label{eq:spaces}
\VV := L^2(0,\ell) \times L^2(0,\ell) 
    \qquad \text{and} \qquad 
\WW :=H^1(0,\ell) \times H^1(0,\ell),
\end{align}  
where $H^1(0,\ell)$ 
denotes the 
standard Sobolev space on the interval $(0,\ell)$. 
%
%
We further introduce the set of admissible states
\begin{align*}
    \DD = \{(\rho,w) \in \VV : (A1)-(A2) \text{ are satisfied}\}.
\end{align*}
Let us recall that classical solutions of  \eqref{eq:sys1}--\eqref{eq:sys3} satisfy $\bu=(\rho,w) \in C^1([0,T];\VV)\cap C^0([0,T];\DD)$ and $\bz(\bu) = (h,m) \in C^0([0,T];\WW)$. 
Then $\C$ and $\R(\bu)$ with $\bu \in \DD$ can be understood as self-adjoint 
and positive semi-definite bounded linear operators mapping from $\VV$ or $\WW$ to the dual spaces $\VV'$ or $\WW'$. 
For $\eps$, $a$ uniformly positive and bounded, $\C$ induces a norm 
\begin{align} \label{eq:normC}
    \|\bu\|_\C^2 =  \|\sqrt{a} \rho\|^2_{L^2(0,\ell)} +  \|\eps w\|^2_{L^2(0,\ell)},
\end{align}
which is equivalent to the standard norm on $\VV = L^2(0,\ell) \times L^2(0,\ell)$. 
Adopting the previous notation, we write $\la \cdot, \cdot\ra$ for the duality products on $\VV' \times \VV$ and $\WW' \times \WW$, respectively, and use $(a,b)=\int_0^\ell a(x) b(x) \, dx$ to denote the scalar product of $L^2(0,\ell)$. 
From the variational definition of the operator $\J$, one can see that
\begin{align} \label{eq:skew}
    \la \J \bw, \widetilde \bw \ra:= 
    & (\dx m, \tilde h) - (h, \dx \tilde m) = -\la \J \widetilde \bw, \bw\ra, 
\end{align}
for all $\bw=(h,m)$, $\widetilde \bw=(\tilde h, \tilde m) \in \WW$, i.e., $\J$ is skew-symmetric on $\WW$. 
The formula $\la \B_\partial \bz(\bu), \bv\ra = - h  r\big|_0^\ell$ with $\bv=(q,r)$ finally shows that the boundary operator $\Bd$ acts on the co-state variables. Note that the required boundary values are well-defined for functions $\bz(\bu) = (h,m)$ and 
$\bv=(q,r) \in  \WW = H^1(0,\ell) \times H^1(0,\ell)$.

In summary, we thus have shown that \eqref{eq:sys1}--\eqref{eq:sys3} can be interpreted as an abstract port-Hamiltonian system \eqref{eq:abs1}--\eqref{eq:abs2}, and under conditions (A1)--(A2) also Assumption~\ref{ass:1} is valid.

\section{Stability analysis and parabolic limit} \label{sec:4}

In the following, we verify the conditions of Assumption~\ref{ass:2} for the gas transport problem \eqref{eq:sys1}--\eqref{eq:sys3}, and then utilize the abstract stability results to prove stability of solutions with respect to perturbations in the model parameters as well as initial and boundary values. As a case of particular interest, we will study  convergence in the parabolic limit $\eps \to 0$. 

\subsection{Perturbed problem}

Let $(\hat \rho,\hat w)$ denote a solution to the perturbed equations
\begin{align}
    a \dtau \hat \rho + \dx (a \hat \rho \hat w) &= 0,  \label{eq:pert1} \\
    \hat \eps^2 \dtau \hat w + \dx ( P'(\hat \rho) + \hat \eps^2 \frac{\hat w^2}{2} + g z)  + \hat \gamma \frac{|\hat w|}{ a\hat \rho} \hat m &= 0, \label{eq:pert2}
\end{align}
which are again assumed to hold for all $0 < x < \ell$ and time $t>0$. 
The terms carrying spatial derivatives are the corresponding co-state variables
\begin{align} \label{eq:pert3}
    \hat h(\hat \rho,\hat v) = \hat \eps^2\frac{\hat w^2}{2} + P'(\hat \rho) +  g z 
\qquad \text{and} \qquad 
    \hat m(\hat \rho,\hat v) = a \hat \rho \hat w,
\end{align}
which are again directly related to the derivatives of the corresponding energy functional $\hat \H(\hat \rho,\hat w) = \int_0^\ell a (\hat \eps^2 \frac{\hat w^2}{2\hat \rho} + P(\hat \rho) + \hat \rho g z ) \, dx$. 
By the some elementary manipulations and the same arguments as employed above, this system can again be written in the abstract form
\begin{align}
    \C \dtau \buh + (\J + \R(\buh)) \bz(\buh) &= \Bd \bz(\buh) + \beh \\
    \bz(\buh) &= \C^{-1} \H'(\buh),
\end{align}
with the functionals $\H(\cdot)$ and $\bz(\cdot)$, and the operators $\C$, $\J$, and $\R(\cdot)$ denoting those for the unperturbed problem with parameters 
$\eps$ and $\gamma$,
%
and residual $\beh=(\beh_1,\beh_2)$ given by
 %
\begin{align} \label{eq:beh}
\beh_1 = 0
\qquad \text{and} \qquad
\beh_2 = (\eps^2- \hat\eps^2) (\dtau \hat w  + \tfrac{1}{2} \dx |\hat w|^2) +  (\gamma - \hat \gamma) |\hat w| \hat w. 
\end{align}
In the following section, we verify the abstract assumptions required for our stability analysis, without explicitly taking into account the special form of $\beh_1$ and $\beh_2$. 
%

\subsection{Verification of conditions (C0)--(C4)} \label{sec:aux}

We always assume in the following that assumptions (A1)--(A2) are valid. Constants arising in the estimates may depend on the bounds of these conditions.
Since we later consider the case $\eps \to 0$, we will make explicit the dependence of the constants on this scaling parameter.

\subsection*{Condition (C0)}

From Taylor's formula, we know that 
\begin{align*}
f(\bu|\hat \bu) 
&= f(\bu) - f(\buh) - f'(\buh) (\bu-\buh) \\
&= \frac{1}{2} \int_0^1 \la (1-s) f''(\buh + s (\bu-\buh)) (\bu-\buh), \bu-\buh) \rangle ds.
\end{align*}
Now let $f(\bu) = a(\eps^2 \rho \frac{|w|^2}{2} + P(\rho))$ denote the integrand of the energy functional defined above, and note that its Hessian is given by
\begin{align*}
    f''(\rho,w) = a \begin{pmatrix}
    P''(\rho)  &  \eps^2 w \\  \eps^2 w & \eps^2 \rho
    \end{pmatrix}.
\end{align*}
By multiplying with $\bv=(x,y)$ from the left and right, one can see that
\begin{align*}
 \frac{1}{a} \la f''(\rho,w) \bv, \bv\ra 
 &= P''(\rho) x^2 + 2 \eps^2 w x y + \eps^2 \rho y^2 \\
 &= P''(\rho) x^2 + 2 \eps w x (\eps y) + \rho (\eps y)^2.
\end{align*}
The second term in the second line can be estimated by 
\begin{align*}
|2 \eps w x (\eps y) | 
&\le 2 \frac{\eps^2 w^2}{\rho} x^2 + \frac{1}{2} \rho (\eps y)^2.
\end{align*}
From condition \eqref{eq:ass1}, one can see that  $2 \frac{\eps^2 w^2}{\rho} \le P''(\rho)/2$, which leads to the lower bound
\begin{align*}
\la f''(\rho,w) \bv, \bv\ra \ge \frac{a}{2} \left( P''(\rho) x^2 + \rho |\eps y|^2 \right), \\
\intertext{and the corresponding upper bound}
\la f''(\rho,w) \bv, \bv\ra \le \frac{3a}{2} \left( P''(\rho) x^2 + \rho | \eps y|^2 \right).
\end{align*}
Using the uniform bounds for $\rho$, $P''(\rho)$, and $a$, we arrive at the following result.
\begin{lemma} \label{lem:C0}
Let (A1)--(A2) hold. Then 
\begin{align*}
    \hat c_0 \|\bu - \buh\|_\C^2 \le \H(\bu|\buh) \le \hat C_0 \|\bu - \buh\|_\C^2 
\end{align*}
with positive constants $\hat c_0,\hat C_0$ only depending on the bounds on $\rho$ and $P''(\rho)$ in (A1)--(A2).
\end{lemma}

\subsection*{Condition (C1)}

Using $\bu=(\rho,w)$ and $\buh=(\hat \rho,\hat w)$, and the definition of $\R(\bu)$ for our particular problem, the left hand side of \eqref{eq:term1} can be expressed as
\begin{align*} 
-\la \R(\bu) \bz(\bu) &- \R(\buh) \bz(\buh),  \bz(\bu) -\bz(\buh)\ra \\
&= - \int_0^\ell   \gamma \, (|w| w - |\hat w| \hat w) \, a \,   (\rho w - \hat \rho \hat w) \, dx =: (*).
\end{align*}
The first term in the integrand can be written as
$|w| w - |\hat w| \hat w 
= \int_0^1 |\hat w+s(w-\hat w)|  ds (w - \hat w)$,
and one can observe that 
$\frac{|w|+|\hat w|}{4} \le \int_0^1 |\hat w + s (w-\hat w)| ds \le \frac{|w|+|\hat w|}{2}$.
The second term in the integrand can be expanded as $\rho w - \hat \rho \hat w = \hat \rho (w - \hat w) + (\rho - \hat \rho) w $, which leads to
\begin{align*}
(|w| w - |\hat w| \hat w) \,  (\rho w - \hat \rho \hat w) 
&\ge \hat \rho \frac{|w|+|\hat w|}{4} |w-\hat w|^2 -  |w-\hat w| \frac{|w+\hat w|}{2} |w| |\rho - \hat \rho| \\
&\ge \hat \rho \frac{|w|+|\hat w|}{8} |w - \hat w|^2 - (|w|+|\hat w|) \frac{|w|^2}{\hat \rho} |\rho - \hat \rho|^2.
\end{align*}
In summary, we thus arrive at the estimate
\begin{align*}
(*) \le - \frac{1}{8} \int_0^\ell \gamma a \hat \rho (|w|+|\hat w|) |w-\hat w|^2 dx + \int_0^\ell \gamma \frac{|w|^2}{\hat \rho} (|w| + |\hat w|) a |\rho - \hat \rho|^2 dx.
\end{align*}
The uniform bounds for $\hat \rho$, $w$, $\hat w$, $\gamma$, $a$ and condition \eqref{eq:norm} then lead to the following result.
\begin{lemma} \label{lem:C1}
Let assumptions (A1)--(A2) be valid. Then 
condition \eqref{eq:term1} holds with
\begin{align} \label{eq:reldiss3}
\D (\bu|\buh) = \frac{1}{16} \int_0^\ell \gamma a \hat \rho (|w| + |\hat w|) (w- \hat w)^2 ds,
\end{align}
where  $\bu=(\rho,w)$ and $\buh=(\hat \rho,\hat w)$, and with constant $\hat C_1=2 \bar \gamma |\bar w|^2 \ubar \rho^{-1} \hat c_0^{-1}$.
\end{lemma}

\begin{remark}
By the elementary fact that $|w|+|\hat w| \ge |w-\hat w|$, one can see that 
\begin{align} \label{eq:reldiss4}
\D (\bu|\buh) \ge c_D \|w-\hat w\|_{L^3}^3,
\end{align}
with positive constant $c_D=\frac{\ubar \gamma \ubar a \ubar \rho}{16 }$. Thus, the relative dissipation $\D(\bu|\buh)$ provides control over the velocity perturbation even in the case $\eps \to 0$, where the velocity contribution to the relative energy $\H(\bu|\buh)$ disappears. This will later be used in our stability analysis. 
\end{remark}

\subsection*{Condition (C2)}

Let us start with noting that 
\begin{align}\label{eq:zzG}
\bz(\bu) - \bz(\buh) - \G(\buh)  (\bu - \buh)
= \begin{pmatrix}
  P'(\rho|\hat \rho) + \eps^2 (w - \hat w)^2/2\\
  a (\rho - \widehat \rho)  (w - \widehat w)
  \end{pmatrix},
\end{align}
which follows directly from the definitions of the gradient $\bz(\bu)$ and the Hessian $\G(\bu)$ for the problem under investigation. By assumption (A1)--(A2), the pressure potential $P$ is smooth and $\rho,\hat \rho$ are uniformly bounded, and consequently
\begin{align*}
    P'(\rho|\hat \rho) \le C a |\rho - \hat \rho|^2,
\end{align*}
with some constant $C$ only depending on the bounds of the coefficients, the density, and the pressure potential. 
The second line in \eqref{eq:zzG} can be estimated by 
\begin{align*}
    a (\rho - \widehat \rho)  (w - \widehat w) \le  \frac{1}{\eps} ( a^2 (\rho - \hat \rho)^2 + \eps^2 (w - \hat w)^2 )
\end{align*}
via Young's inequality.
The left hand side of \eqref{eq:term2} can then be bounded by 
\begin{align*}
  \la \C \dtau \buh, &\bz(\bu) - \bz(\buh) - \G(\buh) (\bu - \buh)\ra \\
  &= \int_0^\ell a \dtau \hat \rho \big(P'(\rho|\hat \rho) + \eps^2 (w-\hat w)^2/2 \big)dx + \int_0^\ell \eps^2 \dtau \hat w a (\rho - \hat \rho) (w - \hat w) dx \\
  &\le (\|a \dtau \hat \rho\|_{L^\infty}  + \|\eps \dtau \hat w\|_{L^\infty}) ((C+a) |\rho - \hat \rho|^2 + \tfrac{3}{2} \eps^2 |w - \hat w|^2). 
\end{align*}
Together with the bounds \eqref{eq:norm}, we thus obtain the following result.

\begin{lemma} \label{lem:C2}
Let (A1)--(A2) hold. Then \eqref{eq:term2} is valid with $\hat C_2 = C (\|\dtau \hat \rho\|_{L^\infty}  + \|\eps \dtau \hat w\|_{L^\infty})$ and constant $C \ge 0$ depending only on the bounds in the assumptions. 
\end{lemma}

\subsection*{Condition (C3)}

We start by expanding 
\begin{align*}
 \la \beh, \bz(\bu) - \bz(\buh)\ra
 &= \int_0^\ell \beh_1 \left(\frac{\eps^2}{2} (|w|^2-|\hat w|^2) + P'(\rho) - P'(\hat \rho) \right) dx \\ 
 & \qquad \qquad \qquad  + \int_0^\ell \beh_2 a (\rho w - \hat \rho \hat w) \, dx = (i)+(ii).
\end{align*}
Using the uniform bounds in assumption (A2) and smoothness of $P(\cdot)$, we obtain 
\begin{align*}
    (i) &\le (\bar w \eps^2 \|w - \hat w\|_{L^2}  + C_P'' \|\rho - \hat \rho\|_{L^2}) \|\beh_1\|_{L^2}  \\
    &\le \eps^2 \|w-\hat w\|^2_{L^2} + \frac{1}{2}\|\sqrt{a}(\rho - \hat \rho)\|^2_{L^2} + C_1 \|\beh_1\|^2_{L^2},
\end{align*}
where we applied H\"older's inequality in the last step.
Condition \eqref{eq:norm} then allows to bound the first two terms by the relative energy.
For the second term, we obtain
\begin{align*}
(ii) &= \int_0^\ell \beh_2 a ( (\rho - \hat \rho) w +\hat \rho (w-\hat w)) dx \\
&\le \bar w \sqrt{\bar a} \|\beh_2\|_{L^2} \|\sqrt{a}(\rho - \hat \rho)\|_{L^2}  +   \bar a \bar \rho \| \beh_2\|_{L^{3/2}} \|(w - \hat w)\|_{L^3} \\
&\le \frac{1}{2}\|\sqrt{a}(\rho - \hat \rho)\|^2  + C_2 \|\beh_2\|_{L^2}^2 + \delta  \|w-\hat w\|_{L^3}^3 + C(\delta) \|\beh_2\|_{L^{3/2}}^{3/2}. 
\end{align*}
Choosing $\delta=c_D$ and using \eqref{eq:reldiss4} allows to bound the third term in the last line by the relative dissipation. 
In summary, we then obtain the following result. 
\begin{lemma} \label{lem:C3}
Let (A1)--(A2) be valid. Then condition \eqref{eq:term3} holds with $\hat C_3 = 1$ and
\begin{align*}
 \P(\beh) = C_1 \|\beh_1\|_{L^2}^2 + C_2 \|\beh_2\|_{L^2}^2 + C_3 \|\beh_2\|_{L^{3/2}}^{3/2},  
\end{align*}
with constants $C_1$, $C_2$, and $C_3$ only depending on the bounds in (A1)--(A2).
\end{lemma}

Using the specific form of the residual given in \eqref{eq:beh}, we may further estimate the perturbation functional $\P(\beh)$ as follows. 
\begin{corollary} \label{cor:C3}
Let the conditions of Lemma~\ref{lem:C3} be valid and $\beh$  be defined as in \eqref{eq:beh}. Then 
\begin{align*}
    \P(\beh) \le \hat C_3' |\eps^2 - \hat \eps^2|^{3/2} + \hat C_3'' |\gamma - \hat \gamma|^{3/2},
\end{align*}
with constants $\hat C_3'$, $\hat C_3''$ only depending on the bounds in assumptions (A1)--(A2) as well as on $\|\dtau \hat w\|_{L^\infty(0,T;L^2)}$ and $\|\dx \hat w\|_{L^\infty(0,T;L^2)}$. 
\end{corollary}

\subsection*{Condition (C4)}


Using the definition of the state and co-state variables, as well as the variational characterization of the boundary operator, we immediately obtain 
\begin{align*}
\la \Bd (\bz(\bu) &- \bz(\buh)), \bz(\bu) - \bz(\buh) \ra  \\
&= - (h(\rho,w) - h(\hat \rho,\hat w)) (m(\rho,w) - m(\hat \rho,\hat w))\big|_0^\ell \\
&\le |h(\rho,w) - h(\hat \rho,\hat w)|_\partial |m(\rho,w) - m(\hat \rho,\hat w)|_\partial.
\end{align*}
with $h(\tilde \rho,\tilde w) = \eps^2 \frac{\tilde w^2}{2} + P'(\tilde \rho)$ and $m(\tilde \rho,\tilde w)=a\tilde \rho \tilde w$ denoting the co-state mappings of the unperturbed problem, and $|a|_\partial = \sqrt{|a(0)|^2 + |a(\ell)|^2}$ denoting the $\ell^2$-scalar product on the space of boundary values. 
We thus obtain
\begin{lemma} \label{lem:C4}
Let assumptions (A1)--(A2) hold. 
Then condition \eqref{eq:term4} is valid with 
\begin{align*}
 \Pd(\bz(\bu)-\bz(\buh)) 
 &= |h(\rho,w) - h(\hat \rho,\hat w)|_\partial |m(\rho,w) - m(\hat \rho,\hat w)|_\partial.
\end{align*}
\end{lemma}

While $h(\rho,w)$ and $m(\rho,w)$ amount to the natural boundary values of the unperturbed problem, the evaluation $h(\hat \rho,\hat w)$ and $m(\hat \rho,\hat v)$ of the unperturbed co-state mappings at the perturbed states does not have a physical meaning. We therefore decompose 
\begin{align*}
  h(\hat \rho,\hat v) &= \hat h(\hat \rho,\hat v) + (h(\hat \rho,\hat v) - \hat h(\hat \rho,\hat v))
\end{align*}
into the natural boundary value $\hat h(\hat \rho,\hat w) = \hat \eps^2 \frac{\hat w^2}{2} + P'(\hat \rho)$ and a corresponding perturbation $h(\hat \rho,\hat v) - \hat h(\hat \rho,\hat v)$ of the state to co-state mapping.
The latter can be estimated by the bounds in assumptions (A1)--(A2), leading to the following result.
\begin{corollary} \label{cor:C4} 
Let the assumptions of Lemma~\ref{lem:C4} hold and $h_\partial=h(\rho,v)|_{\{0,\ell\}}$ and $\hat h_\partial=\hat h(\hat \rho,\hat v)|_{\{0,\ell\}}$ denote the boundary values of the co-state variables of the unperturbed and perturbed system, respectively.
Then  condition \eqref{eq:term4} holds with 
\begin{align*}
 \Pd(\bz(\bu)-\bz(\buh)) 
 &= \hat C_\partial \left( |h_\partial - \hat h_{\partial}| + |\eps^2 - \hat \eps^2| \right)
\end{align*}
with $C_\partial$ depending only on the bounds in (A1)--(A2) and the Lipschitz bounds for $(\hat \rho,\hat w)$. 
\end{corollary}

\subsection{Stability estimate}

Having verified the conditions of our abstract stability analysis, we can now apply Lemma \ref{lem:abs_st} to obtain the following stability estimate.

\begin{theorem} \label{thm:main}
Let assumptions (A1)--(A2) hold and let $(\rho,w)$ and $(\hat \rho, \hat w)$ denote corresponding classical solutions of \eqref{eq:sys1}--\eqref{eq:sys3} for parameters $\eps,\gamma$ and $\hat \eps,\hat \gamma$, respectively. 
Further assume that $(\hat \rho,\hat w)$ is Lipschitz continuous. 
Then 
\begin{align*}
\|\rho(\tau) - \hat \rho(\tau)\|_{L^2(0,\ell)}^2  &+ \eps^2 \|w(\tau) - \hat w(\tau)\|_{L^2(0,\ell)}^2 + \int_0^\tau \| w(s) - \hat w(s)\|_{L^3(0,\ell)}^3 ds\\ 
&\le \hat C e^{\hat c \tau}
\big(\|\rho(0) - \hat \rho(0)\|_{L^2(0,\ell)}^2  + \eps^2 \|w(0) - \hat w(0)\|_{L^2(0,\ell)}^2 \\
& \qquad\qquad\qquad \qquad   + |\gamma - \hat \gamma|^{3/2} + |\eps^2 - \hat \eps^2| + \int_0^\tau |h_\partial(s) - \hat h_\partial(s)|_\partial ds \big),
\end{align*}
where $h_\partial = h(\rho,w)|_{\{0,\ell\}}$ and $\hat h_\partial = \hat h(\hat \rho,\hat w)|_{\{0,\ell\}}$ are the boundary values of the corresponding co-state variables.
Moreover, the constants $\hat c$, $\hat C$ in this estimate only depend on the bounds in assumptions (A1)--(A2) and the Lipschitz bounds for $(\hat \rho,\hat w)$.
\end{theorem}

\begin{remark}
Let us briefly discuss the conditions and conclusions of the theorem: The additional regularity for the reference solution $(\hat \rho,\hat w)$ is required in Lemma~\ref{lem:C2} and Corollary~\ref{cor:C3}.
The reduction in the convergence rate with respect to $\gamma$ comes from the fact, that the error in $w$ is estimated by the friction term rather than the kinetic energy, in order to obtain estimates that are uniform in $\eps$.
The further reduction of the convergence order in $\eps$ is due to effects from the boundary.
With minor modifications of the proofs, one could handle further perturbations, e.g., in the pressure potential $P(\cdot)$, the cross-section $a$, or the height function $z$, and also deal with other boundary conditions.
A careful inspection of the proofs would also allow to relax the smoothness assumptions on $(\rho,w)$ and $(\hat \rho,\hat w)$ to some extent. 
\end{remark}

\subsection{The parabolic limit problem}

We now study convergence of solutions to \eqref{eq:sys1}--\eqref{eq:sys3} in the limit $\eps \to 0$. 
For $\hat \eps=0$ and $\hat \gamma=\gamma$, the resulting limit problem reads 
\begin{align} \label{eq:par1}
    a \dt \hat \rho + \dx (a \hat \rho \hat w) &= 0, \\
    \dx (P'(\hat \rho)) + \gamma |\hat w| \hat w &= 0. \label{eq:par2}
\end{align}
This is a degenerate parabolic problem, whose solvability can be deduced from the results in \cite{Bamberger79,Raviart70}; a detailed analysis can be found in \cite{SchoebelKroehn2020}.
A formal application of Theorem~\ref{thm:main} to this limiting case directly leads to the following result.
\begin{theorem} \label{thm:main2}
Let (A1)--(A2) hold and $(\rho,w)$, $(\hat \rho,\hat w)$ denote classical solutions of \eqref{eq:sys1}--\eqref{eq:sys3} and \eqref{eq:par1}--\eqref{eq:par2}, respectively. 
%
%
Further assume that the initial and boundary values coincide, i.e., $\rho=\hat \rho$ at time $t=0$ and $\eps^2 \frac{w^2}{2} + P'(\rho) = P'(\hat \rho)$ at the boundary $x\in\{0,\ell\}$, and that $(\hat \rho,\hat w)$ is Lipschitz continuous.
Then 
\begin{align*}
\|\rho(\tau) - \hat \rho(\tau)\|_{L^2(0,\ell)}^2  + \int_0^\tau \| w(s) - \hat w(s)\|_{L^3(0,\ell)}^3 ds
&\le \hat C e^{\hat c \tau} \eps^2  ,
\end{align*}
with constants $\hat c,\hat C$ having the same properties as in Theorem~\ref{thm:main}.
\end{theorem}

\begin{remark}
On bounded time intervals, the quadratic norm difference between solutions of the hyperbolic problem \eqref{eq:sys1}--\eqref{eq:sys3} and the parabolic limit problem \eqref{eq:par1}--\eqref{eq:par2} is thus bounded by $O(\eps^2)$.
Let us note that a formal asymptotic analysis would predict a rate $O(\eps^4)$ and such a rate has actually been proven in \cite{LattanzioTzavaras2013} for linear friction and unbounded domains. 
As mentioned in the previous remark, the reduction of the convergence order in our results is due to the nonlinear friction term and the perturbations coming from the boundary.
\end{remark}

\section{Extension to gas networks} \label{sec:5}

As a next step, we now show that, using the underlying abstract framework, the results of the previous section can be extended almost verbatim to gas networks. 
We start with extending the rescaled model
\eqref{eq:sys1}--\eqref{eq:sys6} to gas networks and then discuss the modifications needed in the stability and asymptotic analysis.

\subsection{Network topology}

%
Let $(\V,\E)$ denote a directed and connected finite graph with vertices $v \in \V$ and edges $e \in \E$, which are identified with intervals $(0,\ell^e)$.
We denote by $\E(v)$ the set of edges incident to the vertex $v$, and decompose $\V=\V_0 \cup \V_\partial $ into the sets of interior and boundary vertices, characterized by $\V_0=\{v \in \V: |\E(v)|>1\}$ and $\V_\partial=\{v \in \V: |\E(v)|=1\}$. Here $|\E(v)|$ denotes the cardinality of the set $\E(v)$. 
We further associate to any vertex $v \in \V$ and edge $e \in \E(v)$ a number 
\begin{align*}
n^e(v) = 
\begin{cases}
1 & \text{if } e=(\cdot,v), \\
-1 & \text{if } e=(v,\cdot).
\end{cases}
\end{align*}
The vertex $v$ thus corresponds to the end point $\ell^e$ or the start point $0$ of the interval $(0,\ell^e)$ representing the edge $e$. 

\subsection{Gas transport on networks}

After rescaling as outlined in the introduction, also see Appendix~\ref{sec:6}, the gas transport on every edge of the network is described by 
\begin{align}
    a^e \dtau \rho^e + \dx m^e &= 0 , \qquad e \in \E \label{eq:net1}\\
    \eps^2 \dtau w^e + \dx h^e + \gamma^e |w^e| w^e &= 0, \qquad e \in \E. \label{eq:net2}
\end{align}
By superscript $^e$ we here denote functions or parameters restricted to the edge $e$. The corresponding co-state variables are defined accordingly by 
\begin{align} \label{eq:net3}
h^e &= \eps^2 \frac{|w^e|^2}{2} + P'(\rho^e) + g z^e, \qquad e \in \E, \\
m^e &= a^e \rho^e w^e, \qquad \qquad \qquad \qquad \quad \!  e \in \E. \label{eq:net4}
\end{align}
These equations, which correspond to the conservation of mass and the balance of momentum, describe the flow of gas in the individual pipes. 
As outlined in \cite{Reigstad2014,Egger2018}, the coupling across pipe junctions can be modeled by the following conditions
\begin{align}
\sum_{e \in \E(v)} m^e(v) n^e(v) &= 0, \qquad v \in \V_0, \label{eq:net5}\\
h^e(v) &= h^v, \qquad e \in \E(v), \ v \in \V_0, \label{eq:net6}
\end{align}
which correspond to conservation of mass and continuity of the total specific enthalpy $h$ at pipe junctions. 
Note that $h^v$ thus corresponds to the unique value of the enthalpy at the junction $v \in \V_0$.
A combination of the two conditions allows to show that no energy is produced via flow over junctions \cite{Reigstad2014,Egger2018}. 
Similar to the case of a single pipe, we may again prescribe the enthalpy at the boundary vertices by 
\begin{align} \label{eq:net7}
h^e(v) = h_\partial^v, \qquad v \in \V_\partial.
\end{align}

\subsection{Weak formulation and canonical form}

In a similar manner to Section~\ref{sec:3}, we multiply \eqref{eq:net1}--\eqref{eq:net2} with suitable test functions, integrate over the edges, use integration-by-parts, and then sum over all edges, which immediately leads to the variational equations
\begin{align*}
\sum_{e \in \E} (a^e \dtau \rho^e, q^e)_e + (\dx m^e, q^e)_e &= 0, \\
\sum_{e \in \E} (\eps^2 \dtau w^e, r^e)_e - (h^e, \dx r^e)_e + (\gamma^e \frac{|w^e|}{a^e \rho^e} m^e, r^e)_e &= -\sum_{v \in \V} \sum_{e \in \E(v)} h^e(v) r^e(v) n^e(v),
\end{align*}
which hold for all time independent piecewise regular test functions $q$, $r$, and all $t>0$ of relevance. In the last term, we used the elementary identity 
\begin{align*}
\sum_{e \in \E} \sum_{v \in e} h^e(v) r^e(v) n^e(v) 
&= \sum_{v \in \V} \sum_{e \in \E(v)} h^e(v) r^e(v) n^e(v) 
\end{align*}
to change the order of summation. 
In summary, we see that the weak formulation of \eqref{eq:net1}--\eqref{eq:net4} again amounts to an abstract port-Hamiltonian system  \eqref{eq:abs1}--\eqref{eq:abs2} with operators 
\begin{align*}
    \la \C \bu, \bv\ra &= \sum_e (a^e \rho^r, q^e)_e + (\eps^2 w^e, r^e)_e \\
    \la \R(\bu) \bz, \bv\ra &= \sum_e (\gamma^e \frac{|w^e|}{a^e \rho^e} m^e, r^e)_e, \\
    \la \J \bz, \bv\ra &= \sum_e (\dx m^e, q^e)_e - (h^e, \dx r^e)_e,
\end{align*}
and boundary operator defined by 
\begin{align*}
    \la \B_\partial \bz, \bv\ra &= -\sum_{v \in \V} \sum_{e \in \E(v)} h^e(v) r^e(v) n^e(v).
\end{align*}
Let us note that the coupling and boundary conditions \eqref{eq:net5}-\eqref{eq:net7} have not been incorporated up to this point. 
At this moment, the above variational equations thus describe the gas transport in a collection of separated pipes. 
The corresponding function spaces $\VV$ and $\WW$ for the network are thus simply given by
\begin{align*}
\VV = \prod_{e \in \E} \VV^e 
\qquad \text{and} \qquad \WW = \prod_{e \in \E} \WW^e
\end{align*}  
with $\VV^e = L^2(0,\ell^e) \times L^2(0,\ell^e)$ and 
$\WW^e= H^1(0,\ell^e) \times H^1(0,\ell^e)$
denoting the corresponding spaces for the individual pipes. 

\subsection{Verification of conditions (C0)--(C4)}

By the simple \emph{additive} construction, conditions (C0)-(C3) can be obtained with the same arguments as for a single edge $e \in \E$ and summation over all edges. It thus remains to consider condition (C4) in detail. 

We start with splitting the boundary term via
\begin{align*}
    &\la \Bd (\bz(\bu) - \bz(\buh)), \bz(\bu) - \bz(\buh)\ra \\
    &= -\sum_{v \in \V_0} \sum_{e \in \E(v)} (h^e(\rho^e,w^e) - h^e(\hat \rho^e,\hat w^e)) (m^e(\rho^e,w^e) - m^e(\hat \rho^e,\hat w^e)) n^e \\
    & \quad -\sum_{v \in \V_\partial} \sum_{e \in \E(v)}
    (h^e(\rho^e,w^e) - h^e(\hat \rho^e,\hat w^e)) (m^e(\rho^e,w^e) - m^e(\hat \rho^e,\hat w^e)) n^e
    = (i)+(ii),
\end{align*}
into contributions coming from the junctions $v \in \V_0$ and the boundary vertices $v \in \V_\partial$. 
For the latter, we can use the results for a single pipe derived in Section~\ref{sec:4}, and hence
\begin{align*}
    (ii) \le \hat C_\partial \, ( |h_\partial - \hat h_\partial|_\partial + |\eps^2 - \hat \eps^2| ).
\end{align*}
Note that $|h_\partial|_\partial^2 = \sum_{v \in \V_\partial} |h_\partial(v)|^2$ now is the corresponding $\ell^2$-norm on $\RR^{|\V_\partial|}$.

For the remaining junctions $v \in \V_0$, we proceed as follows:
We now make use of the coupling condition \eqref{eq:net6} and let $h^v$ respectively $\hat h^v$ denote the uniquely determined values of the corresponding co-state variable at the junction $v \in \V_0$. 
Then 
\begin{align*}
    \sum_{e \in \E(v)} (h^e(\rho^e,w^e) &- h^e(\hat \rho^e,\hat w^e)) (m^e(\rho^e,w^e) - m^e(\hat \rho^e,\hat w^e)) n^e \\
    &= \sum_{e \in \E(v)} (h^e(\rho^e,w^e) - h^v) (m^e(\rho^e,w^e) - m^e(\hat \rho^e,\hat w^e)) n^e \\
    &\qquad + \sum_{e \in \E(v)} (h^v - \hat h^v ) (m^e(\rho^e,w^e) - m^e(\hat \rho^e,\hat w^e)) n^e \\
    &\qquad + \sum_{e \in \E(v)} (\hat h^v - h^e(\hat \rho^e,\hat w^e)) (m^e(\rho^e,w^e) - m^e(\hat \rho^e,\hat w^e)) n^e \\
    &= (a) + (b) + (c).
\end{align*}
Due to the coupling condition \eqref{eq:net6}, the first term (a) vanishes. 
Since $h^v$ and $\hat h^v$ are both single valued on $v \in \V_0$ and $m^e(\hat \rho^e,\hat w^e) = \hat m^e(\hat \rho^e,\hat w^e)$, we further obtain 
\begin{align*}
    (b) = (h^v - \hat h^v) \sum_{e \in \E(v)} (m^e(\rho^e,w^e) - \hat m^e(\hat \rho^e,\hat w^e)) = 0,
\end{align*}
where we used the second coupling condition \eqref{eq:net5} for the perturbed and unperturbed problem, respectively. 
Using \eqref{eq:net6} for the perturbed problem, we finally get \begin{align*}
    (c) &= \sum_{e \in \E(v)} (\hat h^e(\hat \rho^e,\hat w^e) - h^e(\hat \rho^e,\hat w^e)) (m^e(\rho^e,w^e) - m^e(\hat \rho^e,\hat w^e)) n^e 
    \le C |\eps^2 - \hat \eps^2|,
\end{align*}
with a constant $C$ depending only on the bounds $\bar a, \bar \rho$ and $\bar w$ in assumption (A1)--(A2).
By summation over all edges, we thus obtain the following result.
\begin{lemma} \label{lem:C3net}
Let (A1)--(A2) hold uniformly for all edges $e \in \E$. Furthermore, let $\bu=(\rho,w)$ and $\hat \bu=(\hat \rho,\hat w)$ denote classical solutions of \eqref{eq:net1}--\eqref{eq:net7} with parameters and data $\eps,\gamma,h_\partial$ and $\hat \eps,\hat \gamma, \hat h_\partial$, respectively, and assume that $(\hat \rho,\hat w)$ is Lipschitz on every edge $e \in \E$. 
Then condition (C4) holds with perturbation functional
\begin{align*}
    P_\partial(\bz(\bu)-\bz(\buh)) = \hat C_\partial (| h_\partial - \hat h_\partial| + |\eps^2 - \hat \eps^2|),
\end{align*}
and $\hat C_\partial$ depending only on the bounds in  (A1)--(A2) and the Lipschitz bounds for $(\hat \rho,\hat w)$.
\end{lemma}

\subsection{Stability and parabolic limit for gas networks}

By the considerations of the previous sections, we immediately see that the assertions of Theorem~\ref{thm:main} and \ref{thm:main2} and also the remarks concerning possible generalizations remain valid verbatim also for gas networks.
For the parabolic limit $\eps \to 0$, we state the corresponding result in detail.
\begin{theorem} \label{thm:main3}
Let assumptions (A1)--(A2) hold and let $(\rho,w)$, $(\hat \rho,\hat w)$ denote classical solutions of \eqref{eq:net1}--\eqref{eq:net7} with $\eps>0$, $\hat \eps=0$, and $\gamma=\hat \gamma$, $h_\partial=\hat h_\partial$. Further assume that $(\hat \rho,\hat w)$ is Lipschitz on every edge $e \in \E$. 
Then 
\begin{align*}
\|\rho(\tau) - \hat \rho(\tau)\|_{L^2(\E)}^2  + \int_0^\tau \| w(s) - \hat w(s)\|_{L^3(\E)}^3 ds
&\le \hat C e^{\hat c \tau} \eps^2  ,
\end{align*}
with constants $\hat c,\hat C$ of the same form as in Theorem~\ref{thm:main} and \ref{thm:main2}. 
\end{theorem}
As similar generalization can also be made for the stability estimate of Theorem~\ref{thm:main}, and also the remarks concerning possible generalizations made for a single pipe generalize almost verbatim to gas networks.

\section{Summary}

In this paper, we studied the stability of classical solutions to hyperbolic balance laws describing the gas transport in pipelines and pipeline networks. 
Our analysis was based on the formulation of a rescaled set of equations as an abstract port-Hamiltonian system involving state and co-state variables. 
Under some general smoothness assumptions on possible solutions, we established quantitative perturbation bounds via relative energy estimates. 
As a particular application, we proved convergence of sufficiently smooth solutions to solutions of the parabolic limit problem in the asymptotic high-friction, long-time, low-Mach limit $\eps \to 0$. 
A key ingredient of our analysis was the proper treatment of boundary conditions, which allowed us to generalize our results 
almost verbatim to gas-networks. 
Natural next steps for future investigation are the structure-preserving discretization and discretization error analysis, which  can most probably be done with similar arguments. 
Another topic of interest might be the further investigation of the parabolic limit problem, which seems to be widely used in the gas network community.

\appendix

\section{Derivation of the model equations} \label{sec:6}

For convenience of the reader, we present a brief derivation of the model equations \eqref{eq:sys1}--\eqref{eq:sys5} starting from the basic equations of gas dynamics, 
and we show that the resulting system can be phrased into the abstract form \eqref{eq:dhs} which will be utilized for our analysis.  

\subsection{Model equations}
The flow of gas through a long pipe, which is identified with the interval $(0,\ell)$ in the sequel, is described by the  balance laws for mass and momentum \cite{BrouwerGasserHerty2011}
\begin{align}
a \dt \rho + \dx (a \rho v) &= 0 \label{eq:sys1a}\\
\dt (a \rho v) + \dx (a \rho v^2 + a p) &= p a_x - \frac{\lambda}{2d} |v| a \rho v - a \rho g z_x. \label{eq:sys2a}
\end{align}
with a quadratic friction law.
Here $\rho$ and $v$ are the density and velocity of the gas, $a=a(x)$ and $d=\sqrt{4 a/\pi}$ are the cross-sectional area and diameter of the pipe, $\lambda=\lambda(x)$ the friction factor, $z=z(x)$ the pipe elevation, and $g$ the gravity constant.
We assume in the following that $a$, $d$, and $z$ are at least Lipschitz continuous in space.
We further assume that the pressure $p=p(\rho)$ is a function of the density $\rho$ alone, and define via 
\begin{align} \label{eq:sys3a}
P(\rho) = \rho \int_1^\rho \frac{p(r)}{r^2} dr,
\end{align}
a corresponding pressure potential. We are thus considering barotropic flow conditions.
For a given state $(\rho,v)$, the total free energy of the system can then be expressed by the Hamiltonian energy functional
\begin{align} \label{eq:energy}
\H(\rho,v) &= \int_0^\ell  a \left(\rho \frac{v^2}{2} + P(\rho) + \rho g z\right) \, dx.
\end{align}
Let us refer to \cite{BrouwerGasserHerty2011} for further details on this model.

\subsection{Transformation to dissipative Hamiltonian form}

From equations \eqref{eq:sys1a}--\eqref{eq:sys1b} and using the elementary identity $\dt (a \rho v) = a \rho \dt v + v a \dt \rho$, one can deduce that 
\begin{align*}
\dt v 
&= \frac{1}{a\rho} \dt (a \rho v) - \frac{v}{a\rho} (a \dt \rho)     \\
&= \frac{1}{a\rho} \left( -\dx (a \rho v^2 + a p(\rho)) + p(\rho) \dx a - \frac{\lambda}{2d} |v| a \rho v - a \rho g z_x + v \dx (a \rho v)\right).
\end{align*}
This expression can be further simplified by using that 
\begin{align*}
v \dx (a \rho v)-\dx (a \rho v^2) 
  &= -a \rho v \dx v,  
   = -\frac a 2 \rho \dx(v^2)\\
-\dx (a p(\rho)) + p(\rho) \dx a &= -a \dx p(\rho), 
\end{align*}
and employing the following elementary relation 
\begin{align*}
\frac{1}{\rho} \dx p(\rho)
  &= \left(\frac{p}{\rho^2} + \frac{p'(\rho)}{\rho} - \frac{p}{\rho^2}\right) \dx \rho 
   = \dx \left(\int_1^\rho \frac{p(s)}{s^2} ds + \frac{p}{\rho}\right)
   =  \dx (P'(\rho))
\end{align*}
between the pressure and the pressure potential. 
Substituting these expressions into the above equation for $\dt v$, we can thus rewrite the system \eqref{eq:sys1a}--\eqref{eq:sys2a} with \eqref{eq:sys3a} into the form
\begin{align} 
a \dt \rho + \dx (a \rho v) &= 0, \label{eq:sys1b}\\
\dt v + \dx \left(\frac{v^2}{2} + P'(\rho) + g z\right) &= -\frac{\lambda}{2d} \frac{|v|}{a\rho} (a\rho v),\label{eq:sys2b}
\end{align}
which is equivalent to \eqref{eq:sys1a}--\eqref{eq:sys2a} for sufficiently smooth solutions.

\subsection*{Co-state variables}

Let $(\rho,v)$ denote a smooth solution of \eqref{eq:sys1b}--\eqref{eq:sys2b}. 
Then  the temporal change of the energy $\H(\rho,v)$ can be expressed by 
\begin{align} \label{eq:ddtH}
\ddt \H(\rho,v) 
&= \int_0^\ell \H_\rho(\rho,v) \, \dt \rho + \H_v(\rho, v) \, \dt v \, dx,
\end{align}
with $\H_\rho(\rho,v)$ and $\H_v(\rho,v)$ denoting the partial derivatives of the \emph{energy-density}, i.e., the integrand in \eqref{eq:energy}, which are here given by
\begin{align} \label{eq:costate}
\H_\rho(\rho,v) = a \left(\frac{v^2}{2} + P'(\rho) + g z\right) 
\qquad \text{and} \qquad  
 \H_v(\rho,v) = a \rho v.
\end{align}
Let us note that these terms also appear in the contributions of \eqref{eq:sys1b}--\eqref{eq:sys2b} carrying spatial derivatives. We can therefore rewrite \eqref{eq:sys1b}--\eqref{eq:sys2b} in more compact form as
\begin{align} 
a \dt \rho + \dx m &= 0 \label{eq:ph1}\\
\dt v + \dx h &= - \frac{\lambda}{2d} \frac{|v|}{a\rho} m \label{eq:ph2} 
\end{align}
with state variables $(\rho,v)$ and \emph{co-state variables} $(h,m)$ given by \eqref{eq:costate}, i.e., 
\begin{align} \label{eq:ph3}
    m = a \rho v \qquad \text{and} \qquad h = \frac{v^2}{2} + P'(\rho) + g z.
\end{align}
Note that the variables $h$ and $m$ have a clear physical meaning, i.e., they denote the \emph{total specific enthalpy} and the \emph{mass flow rate}, respectively. 
The system \eqref{eq:ph1}--\eqref{eq:ph3} yields a closed set of equations which together with appropriate initial and boundary conditions describes the barotropic subsonic flow of gas in a single pipe.

\subsection{Rescaling}

Since the pipes in a gas transport network are usually very long, friction plays a dominant role in practice and the time scales relevant for the transport process are rather large. 
Based on the scaling proposed in \cite{BrouwerGasserHerty2011,Osiadacz87}, we therefore make the substitutions
\begin{align} \label{eq:scaling}
\frac{\lambda}{2d} =\frac{1}{\eps^2} \gamma, 
\qquad 
v = \eps w, 
\qquad \text{and} \qquad
t = \frac{1}{\eps} \tau,
\end{align}
where $\eps$ is some small parameter and $\gamma$, $w$, $\tau$ denote the rescaled friction coefficient, velocity, and time variable, respectively. 
A small parameter $\eps$ thus corresponds to a large-friction, low Mach, and long-time regime. 
With the above substitutions applied, the equations \eqref{eq:ph1}--\eqref{eq:ph3} are transformed into the system \eqref{eq:sys1}--\eqref{eq:sys3} stated in the introduction. 
%

{\footnotesize 
\section*{Acknowledgement}
The authors are grateful for financial support by the German Science Foundation (DFG) via grant TRR~154 (\emph{Mathematical modelling, simulation and optimization using the example of gas networks}), projects~C04 and C05 and the \emph{Center for Computational Engineering} at Technische Universität Darmstadt.
}

\footnotesize 
 \bibliographystyle{abbrv}
 \bibliography{pH}
\end{document}